\documentclass[reqno]{amsproc}

\newtheorem{theo}{Theorem}
\newtheorem{rem}{Remark}

\newtheorem{cor}{Corollary}
\newtheorem{df}{Definition}

\newcommand\eps\varepsilon
\newcommand\ph\varphi
\newcommand\kap\Lambda

\usepackage{enumerate}
\usepackage{graphicx}
\usepackage{float}
\usepackage{yfonts}

\usepackage{amsthm}
\usepackage{amsmath}
\usepackage{amscd}

\begin{document}
\title[ On Some Properties of Measurable Functions]
{On Some Properties of Measurable Functions in Abstract Spaces }

\author[E.  Borisenko]{Evgenii Borisenko\\yevhenii16@yandex.ru\\
 Mechanics and Mathematics Faculty,\\
M. V. Lomonosov Moscow State University\\
Russia, 119899, Moscow,  MGU \\ 
 }
 
\author[ O. Zubelevich]{\\ Oleg Zubelevich$^*$\footnote{$^*$The  corresponding author. The research was funded by a grant from the Russian Science
Foundation (Project No. 19-71-30012)}\\ozubel@yandex.ru\\ 
Steklov Mathematical Institute of Russian Academy of Sciences\\ Mechanics and Mathematics Faculty,\\
M. V. Lomonosov Moscow State University\\
Russia, 119899, Moscow,  MGU 
 }

\date{}


\begin{abstract}In this short note we present several infinite dimensional theorems which  generalize corresponding facts from the finite dimensional  differential inclusions theory.
\end{abstract}

\maketitle
\numberwithin{equation}{section}
\newtheorem{theorem}{Theorem}[section]
\newtheorem{lemma}[theorem]{Lemma}
\newtheorem{definition}{Definition}[section]

\section{Introduction}
In article \cite{fil} a concept of a generalized solution to an ODE with nonsmooth right hand side is introduced.

Namely, the following initial value problem
\begin{equation}\label{sg670}
\dot x=f(t,x),\quad x(0)=\hat x\end{equation}
is considered in some domain $D\subset\mathbb{R}^m,\quad x\in D.$ The function $f$ is 
defined in $[0,T]\times D$ and 
measurable with respect to the standard Lebesgue measure. Moreover it is assumed that there exists  a function $u\in L^1(0,T)$ such that
$$|f(t,x)|\le u(t)$$ holds for almost all $(t,x)\in[0,T]\times D$.

Filippov defines a generalized solution $x(t)$ to  IVP (\ref{sg670}) as an absolutely continuous function in $[0,T]$ such that
$x(0)=\hat x$ and the inclusion
\begin{equation}\label{dfg66}\dot x(t)\in \bigcap_{r>0}\bigcap_{N\in\mathcal N}\mathrm{conv}\, f\Big(t,B_r\big(x(t)\big)\backslash N\Big)\end{equation}
holds for almost all $t\in[0,T].$

Here $\mathrm{conv}\,$ stands for the {\bf closed convex hull}; $B_r(x)$ denotes an open ball of $\mathbb{R}^m$ with the radius $r$ and a center at the point $x$;
and we use $\mathcal N$ to denote a set of all measure-null subsets of $\mathbb{R}^m$.

 Before any  differential equations, this definition states up several questions from real analysis. 
The  main ones are:

1) whether the set in the right-hand side of (\ref{dfg66}) is nonvoid for almost all $t$? 

 2) assume in addition that the function $f$ is continuous in $[0,T]\times D$. Is it true that 
 $$\bigcap_{r>0}\bigcap_{N\in\mathcal N}\mathrm{conv}\, f\Big(t,B_r\big(x(t)\big)\backslash N\Big)=\{f(t,x)\}\,?$$

 Article \cite{fil} answers both questions positively and contains a sketch of the corresponding proofs. 
 Particularly, if the function $f$ is continuous  then the definition of a generalized solution turns into the standard definition of the classical solution to an ODE with a continuous right-hand side.

  Though this sketch is completely correct we believe that it would be useful to provide a detailed proof;
 moreover we get rid of some extra hypotheses and generalize Filippov's assertions to an infinite dimensional case. 
 
 Perhaps such a generalization can be of some interest by itself or for the modern  studies of the  inclusions in infinite dimensional spaces. 

 \section{The Main Theorems}
 Theorems \ref{lemmaInters}, \ref{df660} (see below) generalize Filippov's results which are obtained for the case $X=D,\quad Y=\mathbb{R}^m$.
 \subsection{The Case of Measurable Mapping} 
 Let a set $X$ be equipped with a $\sigma$-algebra and a measure $\mu$.
 
 Let $Y$ stand for a Hausdorff topological space with a countable base. A mapping $f:X\to Y$ is supposed to be measurable with respect to the Borel $\sigma-$ algebra in $Y$.

 \begin{theo}
	\label{lemmaInters}
 	There exists a measurable set $N_0\subset X,\quad \mu(N_0)=0 $ such that  
	\begin{equation}\label{dfbtt7}\bigcap_{\mu(N)=0}{ \overline{f(X\backslash N)}}= \overline{f(X\backslash N_0)}.\end{equation}
	The intersection is taken over all measurable sets of measure zero; the line denotes the closure.
\end{theo}

\begin{cor}
	\label{lemmaConvInters}Assume in addition that  $Y$ is a topological vector space then
		$$\bigcap_{\mu (N)=0}{ \mathrm{conv} \: f(U \backslash N)}= \mathrm{conv} \: f(U \backslash N_0).$$	
\end{cor}
Indeed, the inclusion  $$\bigcap_{\mu (N)=0}{ \mathrm{conv} \: f(U \backslash N)}\subset \mathrm{conv} \: f(U \backslash N_0)$$ is evident:
$N_0$ is one of the sets $N$ from the left side of the formula.
	
	The inverse inclusion is obtained by means of theorem \ref{lemmaInters} with the help of the formula
	$$\mathrm{conv}\: \bigcap_{\alpha}{U_{\alpha}}\subset \bigcap_{\alpha}{\mathrm{conv}\: U_{\alpha}}. $$

	\begin{df}\label{sd12} We shall say that  $y\in Y$ is a bad point if it has an open neighbourhood  $U_y$ such that $\mu(f^{-1}(U_y))=0$.
	We shall call other points of $Y$  good.\end{df}
	
 The notion of a good point has a direct relation to the standard concept of the essential range \cite{rud}.
	
	Generalizing the definition of the essential range to the infinite dimensional case,  we say that the essential range of the function $f$ coincides with a set of the good points:
	$$\mathrm{ess.im}\,f:=\{y\in Y\mid y \,\, \mbox{is good}\}.$$
	
	A finite dimensional version of the following theorem is well known. 
\begin{theo}\label{xsdfg550505}The set in the left hand side of (\ref{dfbtt7}) coincides with the essential range of $f$:
\begin{equation}\label{drg123}\bigcap_{\mu(N)=0}{ \overline{f(X\backslash N)}}=\mathrm{ess.im}\,f.\end{equation} \end{theo}

	\begin{rem} If $X$ is a Hausdorff topological space and $\mu$ is a Borel measure,
	moreover if$$X=Y,\quad f=\mathrm{id}_X$$ then the set of the good points coincides by definition with the support of the measure $\mu$.\end{rem}
	
\begin{theo}\label{sdfg59595}
Let $Q\subset X$ be a measurable set. We use $f_Q$ to denote the following restriction
$$f_Q=f\mid_Q:Q\to Y.$$
Then there is a measurable set $$N_0^Q\subset Q,\quad \mu(N_0^Q)=0,\quad N_0\cap Q\subset N_0^Q$$ such that 
$$\bigcap_{\mu(N)=0}{ \overline{f_Q(Q\backslash N)}}= \overline{f_Q(Q\backslash N_0^Q)}.$$\end{theo}

\subsection{Some Consequence of Theorem \ref{lemmaInters} } Zero-measure subsets of measure spaces have an analog in  topological spaces. This analog is called a first Baire category set.

Let $X$ be a topological space. Recall that a set $A\subset X$ is of the first Baire category if it is contained in a countable union of closed sets and everyone of these sets has empty  interior. Other subsets of $X$ have the second Baire category.

Let $Y$ denote a Hausdorff topological space with a countable base  as above. We do not impose any conditions on a mapping $f:X\to Y$.

Let $\mathcal F\subset 2^X$ stand for a family of the first Baire category sets.

\begin{theo}\label{1r77}There exists a set $E_0\in\mathcal  F$ such that
$$\bigcap_{E\in\mathcal  F}\overline {f(X\backslash E)}=\overline {f(X\backslash E_0)}.$$
\end{theo}

This theorem is a direct consequence from theorem \ref{lemmaInters}. Indeed, to apply  theorem \ref{lemmaInters}  one must point out a $\sigma-$ algebra  and a measure in $X$ such that $f$ becomes a measurable function.

On  a role of the $\sigma-$algebra we take $2^X$ and 
introduce a measure $$\mu:2^X\to [0,\infty]$$ such that
 $$\mu(B)=\left\{\begin{matrix}
0,&\quad \mbox{if}\quad B\in\mathcal F;\\
\infty&\quad\mbox{otherwise.}\\ 
\end{matrix}\right.$$

\subsection{The  Case of the Continuous at a Point Mapping}
\begin{theo}\label{df660}
	Let  $X$ be a topological vector space with the Borel $\sigma$-algebra and a measure $\mu$.
	
	Assume that 
	
	1) for any open nonvoid set $S\subset X$ one has $\mu(S)>0$; 
	
	2) $Y$ is a locally convex Hausdorff space;
	
	3) a mapping  $f:X\rightarrow Y$ is continuous at a point $\tilde{x}\in X$.

	Then it follows that
		$$ \bigcap_{U\in \mathcal{U}} \bigcap_{\mu (N)=0}{ \mathrm{conv} \: f(U \backslash N)}=\{f(\tilde{x})\},$$ where $\mathcal{U}$ is a base of the open neighbourhoods at the point $\tilde{x}$.
\end{theo}
\begin{rem}

1) The function $f$ in theorem \ref{df660} needs not be measurable.

2) If in theorem \ref{df660} we will replace condition 2) with

"$Y$ is a  Hausdorff topological  space"

then the assertion is replaced with
$$ \bigcap_{U\in \mathcal{U}} \bigcap_{\mu (N)=0}{ \overline{ f(U \backslash N)}}=\{f(\tilde{x})\}.$$
If in addition $f$ satisfies the conditions of theorem \ref{lemmaInters} then this formula  implies $f(\tilde x)\in \mathrm{ess.im}\,f.$
\end{rem}
 \section{Proofs}
 \subsection{Proof of Theorem \ref {lemmaInters}}	
	All the points of  $U_y$ are evidently bad provided $y$ is bad. Let  $B$ denote the set of bad points.
		Then we have \begin{equation}\label{xsdfg43}B=\bigcup_{y\in B}{U_y}.\end{equation}
	Here $U_y$ are the neighbourhoods from definition \ref{sd12}.
	
	By the second Lindelof theorem \cite{5}
	one can extract a countable subcovering $$\{U_{y_n}\},\quad\{y_n\}_{n=1}^{\infty}\subset B$$ such that
	$$B=\bigcup_{n\in \mathbb{N}}{U_{y_n}}.$$
	We claim that  $N_0:=f^{-1}(B)$ is the set we are looking for. 
	
	Indeed,  
	\begin{multline*}
		\mu(N_0)=\mu(f^{-1}(B))=\mu\left(f^{-1}\left(\bigcup_{n\in \mathbb{N}}{U_{y_n}}\right)\right) \\
		=\mu\left(\bigcup_{n\in \mathbb{N}}{f^{-1}(U_{y_n})}\right) \leq \sum_{n\in \mathbb{N}}{\mu(f^{-1}(U_{y_n}))}=0.
	\end{multline*}

	The inclusion  $$\bigcap_{\mu (N)=0}{ \overline{f(X\backslash N)}}\subset \overline{f(X\backslash N_0)}$$is trivial: $N_0$ is one of the sets $N$ from the left side of the formula. 
	
	Let us check the backward inclusion "$\supset$".	
	Show that  $\overline{f(X\backslash N_0)}$ does not contain the bad points. 
	
	Indeed, 
	$$f(X\backslash N_0)\cap B=f((X\backslash N_0) \cap f^{-1}(B))=f((X\backslash N_0) \cap N_0)=f(\emptyset)=\emptyset.$$
	Here we employ a   set theory formula \cite{6}:
	$$f(A)\cap B = f\left(A\cap f^{-1}(B)\right).$$
	It follows that $f(X\backslash N_0)\subset Y\backslash B$. By formula (\ref{xsdfg43}) the set $B$ is open, thus $Y\backslash B$ is closed. We consequently obtain \begin{equation}\label{dfg500}\overline{f(X\backslash N_0)}\subset Y\backslash B.\end{equation}
	This means that $\overline{f(X\backslash N_0)}$ consists of the good points only.

	Observe that if $y\in Y\backslash B$ is a good point then for any measure-null set $N\subset X$ we have $$y\in \overline{f(X\backslash N)}.$$
	
	Indeed, take a set  $N,\quad \mu(N)=0$. Let $y$ be a good point.
	This means that any open neighbourhood $V_y$ of the point $y$ is such that 
	$$\mu\big(f^{-1}(V_y)\big)>0.$$
	Thus there exists a point $ x'\in f^{-1}(V_y)\backslash N$. This implies $$ f( x')\in V_y\cap f(X\backslash N).$$ Consequently for any open neighbourhood $V_y$ the set $V_y\cap f(X\backslash N)$ is nonvoid. 
	Therefore we have  $y\in \overline{f(X\backslash N)}$.

	Since  $N$ is an arbitrary measure-null set and $y$ is an arbitrary good point we conclude that the set
	 \begin{equation}\label{dfgooop}\bigcap_{\mu (N)=0}{ \overline{f(X\backslash N)}}\end{equation} contains all the good points and thus
	 $$\overline{f(X\backslash N_0)}\subset\bigcap_{\mu (N)=0}{ \overline{f(X\backslash N)}}.$$
	
	The theorem is proved.
	
	\subsection{Proof of Theorem \ref{xsdfg550505}}
	We have proved that  set (\ref{dfgooop})
	contains all the good points. Due to formula (\ref{dfg500}) we see that  set (\ref{dfgooop}) 
	coincides with the set of the good points.
	This gives (\ref{drg123}).
	
	The theorem is proved.
	
\subsection{Proof of Theorem \ref{sdfg59595}}
Let $B_Q$ stand for a set of bad points of the mapping $f_Q$. It is clear that $B\subset B_Q$. Thus one has
$$N_0\cap Q\subset f^{-1}(B_Q)\cap Q=f^{-1}_Q(B_Q)=N_0^Q.$$
The theorem is proved.
	
 \subsection{Proof of Theorem \ref{df660}}
 First let us check that  for any  $U\in \mathcal{U}$ and for any measurable $N,\quad \mu(N)=0$ it follows that 
	\begin{equation}\label{ryt76}\tilde{x}\in \overline{U \backslash N}.\end{equation}
	Indeed, assume the converse: $\tilde x\notin \overline{U \backslash N}$. Since $ \overline{U \backslash N}$ is a closed set there exists an open set $\tilde U\in \mathcal U$ such that 
	$$\tilde U\cap (\overline{U \backslash N})=\emptyset$$all the more so
	\begin{equation}\label{xfgaa}\tilde U\cap( {U \backslash N})=\emptyset.\end{equation}
	On the other hand we have
	$$\tilde U\cap( {U \backslash N})=(\tilde U\cap U)\backslash(\tilde U\cap N).$$
	Since $\tilde x\in \tilde U\cap U$ the set $\tilde U\cap U$ is open, nonvoid thus $\mu(\tilde U\cap U)>0$;
	and $\mu(\tilde U\cap N)=0.$ Therefore
	$$\mu\big(\tilde U\cap( {U \backslash N})\big)>0.$$
	This contradicts to (\ref{xfgaa}).
	
	Inclusion (\ref{ryt76}) is proved.

Let us prove the inclusion\begin{equation}\label{dfbqq}\bigcap_{U\in \mathcal{U}} \bigcap_{\mu (N)=0}{ \mathrm{conv} \: f(U \backslash N)}\supset \{f(\tilde{x})\}.\end{equation}
  From the properties of continuous functions \cite{6} inclusion (\ref{ryt76}) implies   $f(\tilde{x})\in \overline{f(U \backslash N)}$.
	Now   the following trivial observation
	 $$ \bigcap_{U\in \mathcal{U}} \bigcap_{\mu (N)=0}{ \mathrm{conv} \: f(U \backslash N)}\supset \bigcap_{U\in \mathcal{U}} \bigcap_{\mu (N)=0}{\overline{f(U \backslash N)}}$$
	 implies inclusion (\ref{dfbqq}).

	Let us check the inclusion
	\begin{equation}\label{xdf1111}\bigcap_{U\in \mathcal{U}} \bigcap_{\mu (N)=0}{ \mathrm{conv} \: f(U \backslash N)}\subset \{f(\tilde{x})\}.\end{equation}
	The point $f(\tilde x)$ has a base of convex closed neighbourhoods $\mathcal V$  \cite{7}.

	Thus for any $V\in\mathcal V$ there exists $U\in\mathcal U$ such that $f(U)\subset V$; therefore 
	$$\mathrm{conv}\,f(U\backslash N)\subset V.$$
	 We consequently have 
	$$\bigcap_{U\in \mathcal{U}} \bigcap_{\mu (N)=0}{ \mathrm{conv} \: f(U \backslash N)}\subset V.$$
	
	Since $V$ is an arbitrary element of the base $\mathcal V$ the inclusion (\ref{xdf1111}) is proved.
	
	The theorem is proved.

\end{document}